\documentstyle[leqno,amssymb,
amstex,11pt]{amsart}
\normalsize
\def 	     \n 		{\bf N} 
\newcommand {\e}		{\varepsilon}

\numberwithin{equation}{section}

\newtheorem{th}		          	{Theorem}

\newtheorem{prop}[equation]		{Proposition}

\newtheorem{problem}[equation]			{Problem}

\theoremstyle{definition}

\theoremstyle{remark}
\newtheorem{remark}[equation]		{Remark}

\begin{document}

\title
{Extremal properties of the
first eigenvalue of Schr\"odinger-type operators} 
\author{Lino Notarantonio}
\thanks{The author is a Dov Biegun Postdoctoral Fellow at the Weizmann
Institute of Science}
\address{Department of Theoretical Mathematics,
Weizmann Institute of Science,
Rehovot 76100, Israel}
\email{notaran@@wisdom.weizmann.ac.il}  


\maketitle

\begin{abstract}
Given a separable, locally compact Hausdorff space $X$ and a
positive Radon measure $m(dx)$ on it, we study 
the problem of finding the potential $V(x) \ge 0$ that maximizes 
the first eigenvalue of the Schr\"odinger-type operator
$L+V(x)$; $L$ is the generator of a local Dirichlet form $(a , D[a])$
on $L^2(X, m(dx))$.   
\end{abstract}

\thispagestyle{empty}

\pagestyle{myheadings}
\pagenumbering{arabic}

\section*{Introduction}

\noindent
Let $A > 0$; for 
\[
V \in B_A := \left\{ f : \Big( \int_\Omega |f|^p~ m(dx) \Big)^{1/p}
\le A \right\}, \ \ 1 \le p \le \infty, 
\]
we let $\lambda_1(V)$ denote the first eigenvalue of $L + V(x)$: 
\begin{equation}
\begin{cases}
L u + V u = \lambda_1(V) u, &\ \mbox{in}\ X, \\
u \in D[a].
\end{cases}\label{eqn:my-pbm}
\end{equation}
In this paper we shall be concerned with the following
\\[.0625in]
{\bf Problem:} Determine whether
\begin{itemize}
\item[(1)] the supremum $ \sup \{ \lambda_1(V) :
V \in B_A \} $ is finite;
\item[(2)] there exists $\tilde V \in B_A$ such that
\[
\sup \{ \lambda_1(V) : V \in B_A \} = \lambda_1(\tilde V). 
\]
\end{itemize}

The main assumptions on the local Dirichlet form $(a, D[a])$ are (a1),
(a2), (a3) in \S~1 below. In particular we stress that $(a ,D[a])$
need not be a regular Dirichlet form (according to the terminology in
\cite{fukufuku}).  

\vspace{.0625in}
The paper we mainly refer to, and which inspired the present work, is
\cite{egnell} by H. Egnell (cf. also the references therein for other
contributions to this problem). In \cite{egnell} 
$X$ is a bounded open set in ${\bf R}^d$, 
$d \ge 1$; $m(dx) = k^2 dx$, where $dx$ denotes the Lebesgue
measure, $k \ge 0$ is a measurable function on $X$; the family $B_A$ is
correspondingly defined as above; 
the Dirichlet form 
\[
a[u,u] := \int_X \sum_{i,j=1}^n a_{ij}(x)\frac{ \partial u }{
\partial x_j} 
\frac{ \partial u }{\partial x_i}~dx,
\]
with domain $D[a] = \{ u : a [ u,u ] < +\infty\}$; the matrix 
$(a_{ij})_{i,j}$ is symmetric, coercive, that is, there is a
constant $\Lambda >0$ such that 
\[
a_{ij}(x)  \xi_j \xi_j \ge \Lambda |\xi|^2, \ \mbox {for all }\ x\in 
X,\ \xi \in {\bf R}^n,
\]
and $a_{ij} \in L^1(X)$, $i,j=1 , \ldots, n$.

\vspace{.0625in}
The paper is organized as follows. In the first section we fix the
notation, introduce some definitions and preliminary results regarding
the general theory of Dirichlet forms, and in 
the second section we present our solution to the problem considered.
As in \cite{egnell}, the case $p = \infty$ is trivial (with 
maximal potential $V = A$), while the other two cases $1 < p < \infty$
and $p = 1$ are examined with different approaches. The case $1 < p <
\infty$ is treated with a suitable use of standard methods in the
Calculus of Variations. The remaining case $p = 1$
requires the form $(a, D[a])$ to be strongly local (cf. \S~1) and this
case is examined by the analysis of a  related variational inequality
(cf. Proposition~\ref{vi}); we have thus to generalize some results
from the Theory of Variatonal Inequalities to this framework of
Dirichlet forms (Theorem~\ref{rmk-reg-obs}), which is done in the
Appendix; we point out that the energy measure associated with the
strongly local form $(a, D[a])$ (cf. \S~1) plays an important role in
this generalization. 

\vspace{.125in}
{\sc Acknowledgments.} I wish to express my gratitude to
M.~Solomyak for his continuous interest, many discussions and 
suggestions during the preparation of this paper. I also want to
thank V.A.~Liskevich and M.~R\"ockner for discussions about some
parts of the paper; in particular I'm grateful to V.A.~Liskevich who
made me aware of \cite{almaro}, and to M.~R\"ockner who sent me a copy
of it.  
Also several
discussions with R.~Gulliver on a preliminary version of the paper are
gratefully acknowledged.

\section{Preliminaries \& notation}

\noindent
{\em General notation.\/} $X$ is a locally compact separable Hausdorff
space. For any $E \subset X$, $\overline E$ denotes the closure of $E$
in $X$; also we let $\chi_E(x)$ be the function such that $\chi_E(x) =
1$ if $x \in E$, while $\chi_E(x) = 0$ otherwise in $X$.
$C(X)$ denotes the space of all real-valued  continuous
functions $u$ on $X$. A Borel measure (on $X$) is an additive set
function defined on the $\sigma$-algebra generated by the family of
open sets of $X$; a Radon measure is a Borel measure which is finite
on compact sets and different from zero on open non-empty sets. Unless
otherwise 
specified, all the measures under consideration are non-negative. 
We let $m(dx)$ be a Radon measure whose support is the
whole $X$, consider the real Hilbert space $L^2(X, m(dx))$ and for
$u,v \in L^2(X, m(dx))$ we let $(u,v)$ denote their inner product. We
also consider the Banach space $L^p (X,m(dx))$, $ 1 \le p \le \infty$,
the norm of which is denoted by $\| \cdot \|_p$. 
\\
Given two functions $f,g$ on $X$, we denote by $\max \{ f , g\}(x)$
(respectively $ \min \{ f , g\}(x)$) the pointwise maximum
(respectively minimum) between $f(x)$ and $g(x)$, $x \in X$. 

\vspace{.25in}\noindent
{\em Dirichlet Forms~(\cite{fukufuku}).\/} 
A Dirichlet form $(a, D[a])$ on $L^2(X, m(dx))$ is a symmetric,
non-negative bilinear form 
$a[u,v]$ defined on a dense subspace $D[a]$ of $L^2(X, m(dx))$;
moreover $D[a]$ equipped with the {\em intrinsic norm\/} $\bigl(
a[u,u] + (u,u) \bigr)^{1/2}$ is itself a Hilbert space. Thus the
embedding of Hilbert spaces $D[a] \hookrightarrow L^2(X, m(dx))$ is
continuous.  

\vspace{.125in}\noindent
The following result collects some standard properties of functions in
$D[a]$ which will be used in the following. 

\begin{prop}\label{lattice-prop} Let $(a, D[a])$ be a Dirichlet
form. Then 
\begin{itemize}
\item[(1)] If $u\in D[a]$ then the function $v := \min \{ 1, \max \{
u,0\} \}$ belongs to $D[a]$ and $a [ v, v] \leq a [u,u]$.
\item[(2)] The sequence $ (\max \{ -n, \min \{ u, n \} \} )_n $
converges in $D[a]$ to $u$, as $n \to +\infty$. 
\item[(3)] $a [|u|, |u| ] \le a [ u, u]$, for every $u \in D[a]$. 
\item[(4)] If $u,v \in D[a]$, then $\max\{ u,v \},~
\min \{ u,v \} \in D[a]$.
\end{itemize}
\end{prop}


\vspace{.0625in}\noindent 
The form is {\em local\/} if $a[u,v] = 0$ whenever $u, v \in
D[a]$ have disjoint supports; the form is {\em strongly local\/} if
$a[u,v] = 0 $ whenever $u$ is constant on the support of $v$.

\vspace{.25in}\noindent
We shall consider in the rest of the paper the following
conditions. 
\begin{itemize}
\item[(a1)] The embedding $D[a] \hookrightarrow L^2(X, m(dx))$ is
compact.  
\item[(a2)] ``Urysohn-type Property'': For every compact set $K$ and
each relatively compact open set $G \subset X$, with $K \subset G$,
there exists a function $u \in D[a] $ such that 
\begin{align*}
u & = 1, \ \mbox{on}\ K \\
u & = 0, \ \mbox{on}\ X \setminus G.
\end{align*}
\item[(a3)] $C(X) \cap D[a]$ is a {\em core\/} of the form $(a,
D[a])$, that is, $C(X) \cap D[a]$ is dense in both $D[a]$ (with
respect to the intrinsic norm) and in $C(X)$ (with respect to the
uniform convergence on compact sets). 
\end{itemize}

\vspace{.25in}\noindent
If the form $a[u,v]$ is strongly local, then we write it as follows:
\begin{equation}
a[u,v]=\int_X\mu[u,v](dx), \ \ \ \ u,v \in D[a].
\label{eqn:beurl-deny}
\end{equation}
In (\ref{eqn:beurl-deny}),
$\mu[u,v](dx)$ is a signed Borel measure (the energy measure 
associated with the form $a[u,v]$). The mapping $(u,v) \mapsto
\mu[u,v](dx)$ is a symmetric non-negative bilinear form; moreover we
assume that the energy measure satisfies the following {\em
localization property\/}: If $A \subset X$ is any open set and  $u =
v$ $m$-a.e. on $A$, then 
\[
\chi_A(x) \mu[u,u](dx) = \chi_A(x) \mu[v,v](dx).
\]


As a consequence of this property we have  that 
\begin{equation}
\chi_A(x) \mu[u,v](dx) = 0,
\label{loc-char}
\end{equation} 
whenever $u$ is constant on $A$, for every $v\in D[a]$.

\begin{remark}\label{rmk:almaro}
The formula (\ref{eqn:beurl-deny}) is a particular case of a result by
S.~Albeverio, Z.M.~Ma \& M.~R\"ockner \cite[Theorem 1.1]{almaro}
concerning the representation of Dirichlet forms satisfying (a3) via
an extension of the Beurling-Deny formula. 
\end{remark}

\vspace{.25in}\noindent
{\em The generator of the form.\/} The generator of a Dirichlet form
$a[u,v]$ is the {\em non-negative self-adjoint operator\/} $L$ whose
domain $D[L]$ is dense in the domain $D[a]$ of the form and satisfies 
$$
a[u,v] = (L u, v),
$$
for $u\in D[L]$ and $v\in D[a]$.

\section{Solution of the problem}

\noindent
In this section we present our solution of the problem stated in
the Introduction. For the convenience of the reader, we rewrite the
problem here. Let $A >0$, 
let $p\in [1,\infty]$, and let us consider
\[
B_A := \left\{  f: \| f\|_p \le A \right\}.
\]

\vspace{.125in}\noindent
Let $(a, D[a])$ be a Dirichlet form that satisfies (a1), 
and let $L$ be the generator of the form. 
For $V \in B_A$; 
we denote by $\lambda_1(V)$ the first eigenvalue of the problem
\begin{equation}
\begin{cases}
Lu + V u = \lambda_1(V)  u, \ \mbox{in}\ X, \\
u \in D[a], 
\end{cases}
\label{eqn:scr-pbm}
\end{equation}
that is, $u\in D[a]$ and
\[
a[u,w] + \int_X V u w~m(dx) = \lambda_1(V) \int_X  u w~m(dx),
\]
for every $w \in D[a]$.
\\[.125in]
The problem is the following.

\vspace{.125in}\noindent
\begin{problem}\label{pbm:scr-pbm}
Determine whether:
\begin{itemize}
\item[(1)] the supremum $ \sup \{ \lambda_1(V) :
V \in B_A \} $ is finite;
\item[(2)] there exists $\tilde V \in B_A$ such that
\[
\sup \{ \lambda_1(V) : V \in B_A \} = \lambda_1(\tilde V). 
\]
\end{itemize}
\end{problem}

If  such a potential $\tilde V$ exists, then we call
the pair  $(\tilde u , \tilde V)$ composed by the solution $ \tilde u$
to (\ref{eqn:scr-pbm}) with potential $\tilde V$ 
the {\em extremal pair.\/} 

\vspace{.125in}
Let us associate to (\ref{eqn:scr-pbm}) above the corresponding 
Rayleigh quotient defined by 
\[
R_V(u) := \frac{ \displaystyle{
a [u,u] + \int_X Vu^2 m(dx) }}
{ \displaystyle{
\int_X u^2 m(dx)}}, 
\ u\in D[a], \ u \not= 0. 
\]

Adapting the variational principle to our case, we have that the first
eigenvalue $\lambda_1(V) $ in (\ref{eqn:scr-pbm}) can be determined as
the lower bound of the Rayleigh quotient:
\begin{equation}
\lambda_1 (V) = \inf\{ R_V (u) : u \in D[a], \ u \not= 0\}.
\label{eqn:var-princ}
\end{equation}

\begin{remark}\label{p-infty}
From the variational principle (\ref{eqn:var-princ}) above, we see
that the case $p = \infty$ has the trivial solution $V = A$.
\end{remark}


\vspace{.25in}
Denoting by $q$ the conjugate exponent of $p$ ($p^{-1} + q^{-1} =1$),
let us consider the following functional  
\[
J(u) :=  \frac{ a[u,u] + A \| u\|^2_{2q}}
{\| u \|^2_2}, \ u\in  D[a], \ u \not= 0.
\]
By standard properties of Dirichlet forms
(cf. Proposition~\ref{lattice-prop} above) we have that $J(|u|) \le
J(u)$. Notice that the functional $J(u)$ is such that 
$J(tu) = J(u)$, $t\in R$. 

\vspace{.125in}
Bythe H\"older inequality, we have
\[
R_V(u) \le J(u)
\]
for arbitrary $ u\in  D[a]$, $u \not= 0$. Thus from
(\ref{eqn:var-princ}) we get 
\[
\lambda_1 (V) \le \inf \{ J(u) : u \in D[a], \ u \not= 0\}.
\]
Thus we see that $\sup_{V \in B_A} \lambda_1(V) < +\infty$ whenever
the right-hand side in the above inequality is finite. The next result
shows that this is indeed the case. 

\begin{prop}\label{exis-minim} Let $(a , D[a])$ be a Dirichlet form
that satisfies (a1) and (a2) in \S~1. Then 
the functional $J(u)$ attains its minimum in $D[a]$. 
If moreover $V \ge 0$, then also $R_V(u)$ attains its minimum in
$D[a]$. Furthermore the minimizers for both $J(u)$ and
$R_V(u)$  are non-negative. 
\end{prop}

\begin{pf}
We shall prove the existence of minimizers for $J(u)$, the other case
being analogous. First of all we notice that the functional $J(u)$ is
not identically equal to $+\infty$; this is a consequence of (a2)
and of $m(dx)$ being a Radon measure. Let thus $(u_h)_h$ be a
minimizing sequence normalized so that $\| u_h \|_{2} =1$. 

Let us consider first the case $1 < p < \infty$. 
Then $(u_h)_h$ is bounded in $D[a] \cap L^{2q}(X)$,   
and therefore a subsequence $(u_{h'})_{h'}$ of $(u_h)_h$ will converge
to $u\in D[a] \cap L^{2q}(X)$; as the embedding of
$D[a]$ into $L^2(X , m(dx))$ is compact, then the whole
sequence will converge to $u$ in $L^2(X ,m(dx))$. Now a semicontinuity
argument shows that 
\[
J(u) \le \liminf_{ h\to +\infty} J(u_{h'}),
\]
hence $u$ is a minimizer. As $J(|u|) \le J( u )$
(cf. Proposition~\ref{lattice-prop}), the minimizers are
non-negative. 

Now let us examine the case $ p =1 $. Then the
sequence $(u_h)_h$ is bounded in $D[a] \cap L^\infty (X)$; passing to
a subsequence $(u_{h'})_{h'}$ we may assume that 
$u_{h'} \rightarrow u$ weakly in $D[a]$ (hence strongly in
$L^2(X)$, and in particular $\| u \|_{2} = 1$) and
weak$^*$-$L^\infty(X)$; therefore $ u \in D[a] \cap L^\infty(X)$, and 
\[
J( u ) \le \liminf_{h' \to \infty} J(u_{h'}),
\]
so that $u$ is a minimizer. The last inequality follows from the two
inequalities  
\begin{align*}
a[u , u] & \le \liminf_{h' \to \infty} a[u_{h'} , u_{h'}], \\ 
\| u \|_{\infty} & \le 
\liminf_{h' \to \infty} \| u_{h'} \|_{\infty}.
\end{align*}
Arguing as in the previous case, it is easily seen that the minimizers
are non-negative. 
\end{pf}

The next proposition gives a necessary condition for the existence of a
maximizing potential $\tilde V$.

\begin{prop}\label{ext-couple}
Let $(a,D[a])$ be a local Dirichlet form that satisfies (a1), (a2) in
\S~1; let $\tilde u$ be a minimizer of $J(u)$, and assume that there
exists a function $\tilde V \in B_A$ with $\mbox{supp}( \tilde V)
\subset \mbox{supp}( \tilde u)$ such that 
\[
L \tilde u + \tilde V \tilde u = \lambda \tilde u,
\]
where $\lambda := J ( \tilde u)$ is the minimum value of $J(u)$. 
Then the minimum value of the Rayleigh quotient $R_{\tilde V}(u)$ is
equal to the minimum value of the functional $J(u)$ and $\tilde u$ is
a minimizer for $R_{\tilde V}(u)$: 
\[
\inf  \{ R_{\tilde V} (u) : u \in D[a],\ u \not= 0\} 
 = R_{\tilde V}(\tilde u) = \lambda = J(\tilde u).
\]
\end{prop}

\begin{pf} It is similar to the proof of Lemma 6 in
\cite{egnell} but for the convenience of the reader we present it as
well. Without loss of generality we may assume that $\| \tilde u
\|_{2q} = 1$. Let $v$ be a non-negative minimizer of
$R_{\tilde V}$ and assume that $\lambda' := R_{\tilde V} (v) <
\lambda$. Then 
\begin{align*}
L \tilde u + \tilde  V \tilde u & = \lambda \tilde u, \\
L v + \tilde V v & = \lambda' v.
\end{align*}
(As $v$ is a minimizer of the Rayleigh quotient $R_V(\cdot)$, $v$
satisfies the latter equation, while $\tilde u$ satisfies the former
by hypothesis.) This implies that 
\[
(\lambda - \lambda' ) \int_X \tilde u v m(dx) = 0,
\]
so that $\tilde u v = 0$, $m$-a.e. on $ X $. As $\mbox{supp}(
\tilde V) \subset \mbox{supp}(\tilde u)$ and $\tilde u v = 0$,
$m$-a.e. on $ X $, we get $\tilde V v = 0$ $m$-a.e. in $ X $. This
implies in particular that 
\[
\lambda' = \frac{ a[v,v]}{ \| v \|^2_2}
\]
Let $v_n := (1/n) \min\{ v, n \}$, for $n \in \n$, and let 
\[
\lambda_n := \frac{a[v_n,v_n]}{\| v_n \|^2_2}.
\]
By Proposition~\ref{lattice-prop} $nv_n$ is in 
$D[a]$, and  $nv_n$ converges to $v$ in $D[a]$, hence in
$L^2(X , m(dx))$, as $n \to +\infty$ (recall that $D[a]$ is compactly,
hence continuously, embedded into $L^2(X, m(dx))$); thus
\[
\lim_{h \to + \infty} 
\lambda_n = \lambda' := \frac{ a[v,v]  } { \| v \|^2_2}.
\] 
Notice that, by definition, we get  $v_n
\tilde u =0$ $m$-a.e. on $X$, hence 
by using a well-known result in the theory of Dirichlet forms
(cf. {\em e.g.\/} \cite[Lemma 3.1.4]{fukufuku}) and 
by the local property of the form under consideration we have that
$a[\tilde u, v_n] = 0$. 
\\
Let us examine first the case $ p\in (1,\infty)$. 
Choose $n$ large so that $\lambda_n < \lambda$ and consider
$J(\tilde u + \e v_n)$; we have
\[
J(\tilde u + \e v_n)  = \lambda - \e^2 (\lambda - \lambda_n)
\frac{\int_X v_n^2 m(dx) }{ \int_X \tilde u^2 m(dx)} +
o(\e^2) < \lambda,
\]
for $\e > 0$ small enough, hence a contradiction. As in general $
\lambda \ge \lambda'$, we have thus $\lambda = \lambda'$, hence the
result is proved for $ p\in (1,\infty)$.
\\
Let us consider the remaining case $p = 1$. Using $ \| \tilde u +
v_n \|_{\infty} = \| \tilde u \|_{\infty} = 1$, we
get that
\[
J( \tilde u + v_n)  = \frac{ \lambda \| \tilde u \|^2_{2}
+ \lambda_n \| v_n \|^2_{2}   }{  \| \tilde u
\|^2_{2} +  \| v_n \|^2_{2} } < \lambda,
\]
if $n$ is large enough. This is a contradiction, hence $\lambda' =
\lambda$ also in the case when $p = 1$ and the proof is thus
concluded.
\end{pf}


\subsection{The case $1 < p < \infty$} Now we are in a position to
state and prove the existence of the extremal pair for this case.

\begin{th}\label{main-thm} Let $(a , D[a])$ be a local Dirichlet form
that satisfies (a1), (a2) in \S~1, and let $L$ be the associated
non-negative self-adjoint operator. Let $1 < p < \infty$ and let $q =
\frac{p}{p-1}$ be its conjugate exponent. Then there exists an
extremal pair $(\tilde u,\tilde V)$ that solves
Problem~(\ref{pbm:scr-pbm}); the potential $\tilde V$ is the unique
maximizer of the first eigenvalue of 
Problem~(\ref{pbm:scr-pbm}) in $ \{ f \in L^p(X , m(dx)) : \| f \|_p
\le A \}$; the function $\tilde u$ is a non-negative minimizer of the
functional $J(\cdot)$ and is also the first eigenfunction of
(\ref{eqn:scr-pbm}):
\[
\begin{cases}
L \tilde u + \tilde V \tilde u = \lambda_1(\tilde V) \tilde u, \
\mbox{in}\ X, \\ 
u \in D[a],
\end{cases}
\]
where 
\[
\tilde V = \Big( A \| \tilde u^2 \|_{q}^{1-q}\Big) \tilde
u^{2(q-1)},
\]
and $\lambda_1 (\tilde V)$ is the the maximal first eigenvalue with
\[
 \lambda_1 (\tilde V) = J(\tilde u) = R_{\tilde V} (\tilde u).
\]
\end{th}

\begin{pf}
The functional $J(\cdot)$ is Gateaux-differentiable; for $\phi \in
D[a] \cap L^{2q}(X)$ ($\not=  
\emptyset$, by (a2), and being $m(dx)$ a Radon measure) we have
\begin{align*}
J'_\phi(u) = \frac{2}{\| u \|^2_{2}} \Big(
a[u,\phi]  & + A \| u^2 \|^{1-q}_{q} \int_X |u|^{2(q-1)}
u\phi~m(dx) \\
& - J(u) \int_X u\phi~m(dx) \Big).
\end{align*}
By Proposition~\ref{exis-minim} $J(\cdot)$ has (non-trivial) non-negative
minimizers; thus a minimizer $\tilde u \ge 0$ of $J(\cdot)$ solves
the equation 
\[
\begin{cases}
L \tilde u + \tilde V \tilde u = \lambda_1(\tilde V) \tilde u,\\ 
\tilde u \in D[a],
\end{cases}
\]
that is,
\[
\begin{cases}
\displaystyle 
a[\tilde u , w ] + \int_X \tilde V \tilde u w~m(dx) =  
\lambda_1 (\tilde V) 
\int_X \tilde u w~m(dx),\ \mbox{for every}\ w 
\in D[a]\\ 
\tilde u \in D[a],
\end{cases}
\]
where 
\[
\lambda_1(\tilde V):= J(\tilde u), \ \mbox{and}\  \tilde V := \Big( A
\| \tilde u^2 \|^{1-q}_{q}\Big)  \tilde |u|^{2(q-1)}.
\]
A direct computation shows that $\| \tilde V \|_{p} = A$, hence
$\tilde V \in B_A$, and, by its definition, $\mbox{supp}(\tilde V)
\subset \mbox{supp}(\tilde u)$; thus by Proposition~\ref{ext-couple}
$(\tilde u, \tilde V)$ is the extremal couple and $\lambda_1 (\tilde
V) = J(\tilde u)$ is the maximal first eigenvalue. Notice that $\tilde
u$ is the first eigenfunction corresponding to the eigenvalue
$\lambda_1 (\tilde V)$. As for the uniqueness of the maximizing
potential, it is proven similarly as in \cite[Theorem 16]{egnell}. 
\end{pf}

In the same assumptions and notation of the above theorem we have the
following result.  
\begin{prop}\label{l-infty}
The extremal pair $(\tilde u, \tilde V)$ satisfies 
\begin{align*}
\| \tilde u \|_{ \infty} & \le \left( 
\frac{ \lambda_1 (\tilde V) } {A}  
\right)^{ \frac{p-1}{2}}\| \tilde u \|_{2q} < +\infty, \\
\| \tilde V \|_{\infty} & \le \lambda_1 (\tilde V),
\end{align*}
\end{prop}

\begin{pf}
As $\tilde u$ is a minimizer of $J(\cdot)$, by
Proposition~\ref{exis-minim} $\tilde u \ge 0$, and without loss of
generality we may as well assume that $\| \tilde u \|_{2q} =1$;
thus $\| \tilde u^2 \|_{q}^{1-q} =1$. Let
$c := \displaystyle{ \left( \frac{\lambda_1 (\tilde V)}{A}
\right)^{\frac{p-1}{2}}}$ and 
define $\xi := \tilde u -\min \{ \tilde u, c \}$. Notice that $\xi \ge
0$, $\xi \in D[a]$ (by Proposition~\ref{lattice-prop}) and 
\[
a[\xi,\xi] = a[\xi, \tilde u] = \int_X \big( \lambda_1 (\tilde V)
- A|\tilde u|^{2(q-1)} \big) \tilde u \xi~ m(dx);
\]
observe that, with our choice of $c$,  the integrand is negative when
$\xi >0$; thus $\xi = 0$ and this gives $ \tilde u \le c$. The estimate on
$\tilde V$ follows with a  direct computation by using the estimate on
$\tilde u$. 
\end{pf}

\subsection{The case $ p = 1 $} 
Let us consider 
\begin{equation}
K := \{ v\in D[a] : | v |  \le 1\ \mbox{$m$-a.e. on }\ X
\}. 
\label{eqn:convex}
\end{equation}
Using the properties of the Dirichlet form $a[u,v]$ in
Proposition~\ref{lattice-prop} it is not
difficult to see that $K$ is a non-empty, closed, 
convex set in $D[a]$. 
\\
Let us also consider the functional  
\[
T(v) := \frac{ a[v , v] + A  }{ \| v
\|^2_{2}}, \ v\in K.
\]
Let $\tilde u$ be a minimizer of $J(u)$; as $J(t u) = J(u)$, $t \in
R$, we can assume that $ \| \tilde u \|_{\infty} =1$; thus $\tilde
u$ is also a minimizer of $T(\cdot)$ and  
$J(\tilde u) = T(\tilde u)$. 
\\[.0625in]
We have the following result.

\begin{prop}\label{vi}
Let $(a, D[a])$ be a local Dirichlet form on $L^2(X, m(dx))$ 
that satisfies (a1), (a2) and (a3) in \S~1.
Then $\tilde u$ is a solution of the variational inequality 
\begin{equation}
\begin{cases}
\displaystyle{
a[v , v - u] \ge J(\tilde u) \int_X
u(v- u)~m(dx), \ \forall v \in K, }\\
u \in K.
\end{cases}\label{eqn:vi}
\end{equation}
\end{prop}

\begin{pf}
Let $t \in (0,1)$, $v \in K$; then $\tilde u + t (v - \tilde u) \in
K$, $J(\tilde u) = T(\tilde u)$ and 
\[
T(\tilde u) \le T( \tilde u + t (v - \tilde u)).
\]
A direct computation, similarly as in the proof of \cite[Proposition
12]{egnell}, shows that $\tilde u$ satisfies (\ref{eqn:vi}). 
\end{pf}

\vspace{.125in}\noindent
Now we are in a position to state and prove the main result for the
case $p = 1$.

\begin{th}\label{main-thm-1}
Let $(a , D[a])$ be a strongly local Dirichlet form 
that satisfies (a1), (a2), (a3) in \S~1. Then there exists an extremal
pair $(\tilde u, \tilde V)$ that solves Problem~\ref{pbm:scr-pbm} and
has the  following properties:
\begin{itemize}
\item[(i)] $\tilde u$ is a minimizer of $J(\cdot)$.
\item[(ii)] $\tilde u \ge 0$, $\| \tilde u \|_{\infty} = 1$.
\item[(iii)] Let $I := \{ x \in X : \tilde u(x) = 1 \}$; 
then $m(I) > 0$. 
\item[(iv)] $\tilde u$ is the first
eigenfunction of Problem~\ref{eqn:scr-pbm}, $\tilde V =
\displaystyle{ \frac{A}{m(I)}\chi_I}$ and the maximal first eigenvalue
$\lambda_1 (\tilde V) = 
\displaystyle{ \frac{A}{m(I)}}$.
\item[(v)]  $R_{\tilde V} (\tilde u) = J(\tilde u)$.
\item[(vi)] The potential $\tilde V$ is the unique maximizer of the
first eigenvalue of Problem~\ref{eqn:scr-pbm}.
\end{itemize}
\end{th}

\begin{pf}
By Proposition~\ref{exis-minim} the functional $J(u)$ attains its
minimum in $D[a]\cap L^\infty(X)$, and its minimizers
are non-negative. Let $\tilde u$ be a minimizer of $J(u)$ and without
loss of generality we may assume that $ \| \tilde u \|_{\infty} =1$
(recall that $J(t \cdot ) = J( \cdot )$, for $t \in R$) so that $0 \le
\tilde u \le 1$. By Proposition~\ref{vi} $\tilde u$ is a 
solution of the variational inequality (\ref{vi}). Letting $\lambda =
J(\tilde u)$ and considering the ``obstacle'' equal to the constant
function $\psi =1$,
by Theorem~\ref{rmk-reg-obs} in the Appendix we have that $L
\tilde u  = \lambda \tilde u$, on $X \setminus I$, and $L \tilde u =
0$ 
on $ I $, so that 
$ L \tilde u + 
\chi_I(x) \lambda \tilde u 
= \lambda \tilde u$, that is, 
\[
a[\tilde u, v] + \lambda \int_X \chi_I(x)
\tilde u(x) v(x)~m(dx) 
= \lambda \int_X \tilde u(x) v(x)~m(dx), 
\]
for every $v\in D[a]$; in particular for $v= \tilde u$, and
recalling that $\tilde u = 1$ on $I$, we have
\[
a[\tilde u , \tilde u] + \lambda m(I) 
= \lambda \int_X |\tilde u|^2m(dx).
\]
If $m(I) = 0$ then from the relation above we get that $\lambda$ is
the first eigenvalue of the problem (\ref{eqn:scr-pbm}) with $V=0$, 
$\tilde u$ being the corresponding eigenfunction; thus from the
variational principle we have
\begin{equation}
\frac{ a[ \tilde u , \tilde u ] } { \| \tilde u \|^2_2 } = \lambda;
\label{eqn:contr}
\end{equation}
as $\tilde u$ is also a minimizer for $J(u)$ with $\| \tilde u
\|_\infty = 1$, $\tilde u$ is also a minimizer for $T(u)$ and 
\[
\lambda = T( \tilde u) = \frac{ a[ \tilde u , \tilde u ] + A } { \|
\tilde u \|^2_2 }
\]
but from (\ref{eqn:contr}) we get a contradiction, since $A >0$. 
Therefore $m(I) >0$  and we have that 
\[
\lambda = \frac{ a[\tilde u , \tilde u] + \lambda
m(I)  }{ \| \tilde u\|_{2}^2 } = J(\tilde u),
\]
and this implies $A = \lambda m(I)$. Therefore if we define $\tilde V
:= \frac{A}{m(I)} \chi_I $, then by Proposition~\ref{ext-couple} we
have that  $(\tilde u, \tilde V)$ is an extremal pair, $\lambda_1
(\tilde V):= 
A/m(I)$ is the extremal eigenvalue and $R_{\tilde V} (\tilde u ) =
J(\tilde u)$; 
this proves the first five statements in the theorem. As for the
uniqueness of $\tilde V$, we can argue similarly as in \cite[Theorem
16]{egnell} and conclude the proof.
\end{pf}

\section{Appendix}
In this section we deal with a strongly local Dirichlet form $(a,
D[a])$ that 
satisfies (a1), (a2), (a3) in \S~1. With these assumptions we can,
similarly as in \cite[Chapter 3]{fukufuku}, introduce the notions of
capacity (associated with $(a, D[a])$) and quasi-continuity; in
particular we can associate to each $u \in D[a]$ a
sequence of closed sets $(F_k)_k$ (a ``nest'') such that the union
$\bigcup_k F_k$ is equal to $X$ (with the exception perhaps of a set
of capacity zero) and the restriction of $u$ to $F_k$ is continuous on
$F_k$, $k \in N$ (cf. Theorems 3.1.2, 3.1.3 in \cite{fukufuku}).

\vspace{.125in}
Denote by $\langle \cdot, \cdot \rangle$ the duality pairing between
$D[a]$ and its dual $\bigl( D[a] \bigr)'$. As the form $ a[\cdot,
\cdot ]  $ continuous on $D[a]$, we have that $L \phi$, defined
on $D[a]$ by 
\[
L \phi : v \mapsto a [\phi, v ],
\]
is well-defined as an element in $\bigl( D[a] \bigr)'$ and $ \langle L
\phi , v \rangle = a[\phi, v] $.

\vspace{.125in}
Let $\psi: X \longrightarrow R$ be a quasi-continuous function
and consider
\[
K_\psi:= \{ v \in D[a]: v \ge \psi \ 
\mbox{$m$-a.e. on}\ X \}; 
\]
$K_\psi$ is a closed convex set which we {\em assume\/} to be
non-empty. 
Let us consider the following obstacle problem: 
Given $f \in ( D[a] )' $, find 
\begin{equation}
\begin{cases}
u \in K_\psi, \\
a[u , v - u] \ge \langle f, v-u \rangle , 
\ \forall v\in K_\psi, 
\end{cases}
\label{eqn:obs-pbm}
\end{equation}


\begin{th}\label{rmk-reg-obs}
Assume that there exists a unique solution $\tilde u$ to the obstacle
problem (\ref{eqn:obs-pbm}), and let
$I := \{  x \in X : \tilde u(x) = \psi (x)\}$.
Then 
\[
L \tilde u = f, \ \mbox{on}\ X\setminus I.
\]
Furthermore, if the obstacle function $\psi $ is equal to a constant
function, then
\[
L \tilde u = 0, \ \mbox{on}\ I,
\]
{\em i.e.,\/} $\displaystyle{\int_X \chi_I(x) \mu[\tilde u,v](dx) =
0}$, for every $v \in D[a]$. 
\end{th}

\begin{pf} Adapting some arguments in \cite{ks} (cf. in particular
Definition 6.7 in \cite[Chapter II]{ks}), it can be shown that the set
$X \setminus I$ is open; thus for $x_o \in X \setminus I$, there are
two neighborhoods $U$, $G$ of $x_o$ such that $U \subset \overline U
\subset G \subset X \setminus I$, and without loss of generality we
can assume that $G$ is a relatively compact open set. By (a2), with $K
= \overline U$, there exists a function $\phi $ contained in the
domain of the form such that $ \tilde u > \psi + \phi $; moreover for
any $\zeta \in D[a]$ with support in $U$ there is $\e >0$ such that
\[
\tilde u + \varepsilon \zeta \ge \psi + \frac{1}{2}\phi.
\]
Thus $v = \tilde u + \varepsilon \zeta \in K_\psi$; substituting this
$v$ in (\ref{eqn:obs-pbm}) and dividing by $\e$ we get 
\[
\int_U \mu[\tilde u , \zeta ](dx) \ge \langle f, \zeta \rangle,
\]
for every $\zeta \in D[a]$, with support in $  U $. We can
argue similarly with $ v = \tilde u - \e \zeta$ and get 
\[
\int_U \mu[\tilde u , \zeta ](dx) \ge \langle f, \zeta \rangle;
\]
hence $L \tilde u = f$, in $ X \setminus I$. 
\\
Now assume that the obstacle $\psi$ is a constant function, and
without affecting the generality of the argument that follows we can
assume that $\psi =0$; moreover we can also assume that $I$ is
contained in some relatively compact open 
set $\Omega \subset X$. Let $(F_k)_k$ be the nest associated with
$\tilde u$. Thus, except perhaps for a set of 
arbitrarily small capacity, we can assume that the function $\tilde u$ is
continuous on $\Omega'$ with $\Omega' \subset \overline{\Omega'}
\subset \Omega$, hence uniformly continuous on $\overline
{\Omega'}$. Due to the uniform continuity of $\tilde u$ on $\overline
\Omega'$, for every $\e >0$ we can find an open neighborhood $U_\e$ of
$I \cap \Omega'$ such that $U_\e \subset \{ \tilde u \le \e\}$.
By the Urysohn-type property
(a2), there exists $w_\e \in D[a]$ with compact support in $\Omega'$
such that $w_\e = \e$ on ($\{  \tilde u  \le \e\}$, hence on) $U_\e$;
we define $u_\e := \max \{ \tilde u, w_\e\}$ so that $u_\e 
= \e$ on $ U_\e $, $u_\e \in D[a]$ (cf. Proposition~\ref{lattice-prop})
and $u_\e$ converges to $\tilde u$ in $D[a]$. By the local character
of the energy measure (cf. (\ref{loc-char})) the restriction of the
energy measure to $U_\e$, $\chi_{U_\e} (x)\mu[u_\e,v](dx)$, is equal
to zero, for every $\e >0$. Letting $\e \to 0$ we can conclude the
proof. 
\end{pf}

\end{document}